\documentclass[10pt]{amsart}
\usepackage{amssymb,amsfonts,lineno,amsthm,amscd,stmaryrd}
\usepackage[leqno]{amsmath}

\newtheorem{thm}{Theorem}

\newtheorem{lem}[thm]{Lemma}

\newtheorem*{jensen}{Jensen's Theorem}
\theoremstyle{definition}
\newtheorem{definition}[thm]{Definition}
\newtheorem*{furtherremarks}{Further remarks}

\newcommand{\EQ}[1]{\eqref{eq:#1}}
\newcommand{\LEM}[1]{Lemma~\ref{lem:#1}}

\newcommand{\THM}[1]{Theorem~\ref{thm:#1}}

\newcommand{\R}{\ensuremath{\mathbb{R}}}
\newcommand{\ep}{\varepsilon}
\newcommand{\mmbox}[1]{\quad \mbox{#1} \quad}

\begin{document}

\title[An easy proof of the uniqueness of infinity harmonic functions]{An easy proof of Jensen's theorem on the uniqueness of infinity harmonic functions}
\author{Scott N. Armstrong}
\address{Department of Mathematics,
University of California, Berkeley, CA 94720.}
\email{sarm@math.berkeley.edu}
\author{Charles K. Smart}
\address{Department of Mathematics,
University of California, Berkeley, CA 94720.}
\email{smart@math.berkeley.edu}

\date{\today}


\begin{abstract}
We present a new, easy, and elementary proof of Jensen's Theorem on the uniqueness of infinity harmonic functions. The idea is to pass to a finite difference equation by taking maximums and minimums over small balls.
\end{abstract}

\maketitle


In this short article, we present a new proof of the famous result of R. Jensen \cite{Jensen:1993}, which establishes the uniqueness of viscosity solutions of the infinity Laplace equation
\begin{equation}\label{eq:infinity-laplace}
-\Delta_\infty u : = - \sum_{i,j=1}^n u_{x_ix_j}u_{x_i} u_{x_j} = 0,
\end{equation}
in a bounded domain $\Omega \subseteq \R^n$, subject to a given Dirichlet boundary condition. Our argument is elementary, and other than our use of a well-known equivalence (\THM{CEG}, below), our presentation is self-contained. In contrast, previous proofs employ either intricate viscosity solution methods as well as a deep result of Aleksandrov \cite{Jensen:1993,Barles:2001,Aronsson:2004,Crandall:2007}, or follow a less direct path using probabilistic arguments \cite{Peres:2009}.

\begin{definition}
A viscosity subsolution (supersolution, solution) of \EQ{infinity-laplace} is called \emph{infinity subharmonic (superharmonic, harmonic)}.
\end{definition}

We refer to the survey articles \cite{Aronsson:2004,Crandall:2008} for an introduction to the infinity Laplace equation, as well as the definition of viscosity solution. In this article we do not apply the definition of viscosity solution directly. Instead, we use the notion of \emph{comparisons with cones}.

\begin{definition}
A \emph{cone function} with \emph{vertex} $x_0 \in \R^n$ is a function of the form $\varphi(x) = a + b|x-x_0|$, where $a,b\in \R$. A function $u \in C(\Omega)$ is said to \emph{enjoy comparisons with cones from above} if it possesses the following property: for every open $V\subseteq \R^n$ for which $\bar V \subseteq \Omega$, and every cone function $\varphi$ with vertex $x_0 \in \R^n \setminus V$,
\begin{equation*}
u \leq \varphi \  \ \mbox{on} \ \partial V \quad \mbox{implies} \quad u\leq \varphi \ \ \mbox{in} \ V.
\end{equation*}
We say that $u$ \emph{enjoys comparisons with cones from below} if $-u$ enjoys comparisons with cones from above. Finally, $u$ \emph{enjoys comparisons with cones} if it enjoys comparisons with cones from above and below.
\end{definition}

The following result of Crandall, Evans and Gariepy \cite{Crandall:2001} establishes the equivalence between the notions of infinity harmonic functions and functions enjoying comparisons with cones.

\begin{thm}[\cite{Crandall:2001}]\label{thm:CEG}
A function $u \in C(\Omega)$ is infinity subharmonic (superharmonic) if and only if $u$ enjoys comparisons with cones from above (below).
\end{thm}

A less general version of \THM{CEG} was also obtained by Jensen \cite[Lemma 3.1]{Jensen:1993}. In addition to \cite{Crandall:2001}, elementary proofs of this result can be found in \cite{Aronsson:2004,Crandall:2008}.

We now introduce some notation. For $\ep > 0$ and $x\in \R^n$, we write $B(x,\ep) := \{y\in \R^n : |x-y| < \ep \}$. Let $\Omega_\ep$ denote the set of points $x\in\Omega$ for which $\bar B(x,\ep) \subseteq \Omega$. If $u \in C(\Omega)$ and $x \in \Omega_\ep$, then we denote
\begin{equation*}
S^+_\ep u(x) := \max_{y \in \bar B(x,\ep)} \frac{u(y) - u(x)}{\ep} \mmbox{and} S^-_\ep u(x) := \max_{y \in \bar B(x,\ep)} \frac{u(x) - u(y)}{\ep}.
\end{equation*}


The first step in our proof of Jensen's theorem is the following comparison lemma for a finite difference equation. We adapt an argument due to Le Gruyer \cite{LeGruyer:2007}, who proved a similar lemma for a difference equation on a finite graph.

\begin{lem}
\label{lem:finite}
Assume that $u, v \in C(\Omega) \cap L^\infty(\Omega)$ satisfy the inequalities
\begin{equation}
\label{eq:finite-hyp}
S^-_\ep u(x) - S^+_\ep u(x) \leq 0 \leq S^-_\ep v(x) - S^+_\ep v(x) \quad \mbox{for every} \ \ x\in\Omega_\ep.
\end{equation}
Then
\begin{equation*}
\sup_{\Omega} (u - v) = \sup_{\Omega \setminus \Omega_\ep} (u - v).
\end{equation*}
\end{lem}

\begin{proof}
Arguing by contradiction, we suppose that
\begin{equation*}
\sup_\Omega (u - v) > \sup_{\Omega \setminus \Omega_\ep} (u - v).
\end{equation*}
The set $E := \{ x \in \Omega : (u-v)(x) = \sup_\Omega (u-v) \}$ is nonempty, closed and contained in $\Omega_\ep$. Define $F := \{ x \in E : u(x) = \max_E u \}$, and select a point $x_0 \in \partial F$. Since $u-v$ attains its maximum at $x_0$, we see that
\begin{equation}\label{eq:finite-max-point}
S^-_\ep v(x_0) \leq S^-_\ep u(x_0).
\end{equation}

Consider the case $S^+_\ep u(x_0) = 0$. From the first inequality in \EQ{finite-hyp} we also have $S^-_\ep u(x_0) = 0$. From \EQ{finite-max-point}, we deduce that $S^-_\ep v(x_0) = 0$. Now the second inequality of \EQ{finite-hyp} implies $S^+_\ep v(x_0) = 0$. In particular, $u$ and $v$ are constant in $\bar B(x_0, \ep)$, contradicting our assumption that $x_0 \in \partial F$.

It remains to examine the case $S^+_\ep u(x_0) > 0$. Select $z \in \bar B(x_0,\ep)$ such that $\ep S^+_\ep u(x_0) = u(z) - u(x_0)$. Since $u(z) > u(x_0)$ and $x_0 \in F$, we see that $z\not\in E$. From this we deduce that
\begin{equation}\label{eq:temp}
\ep S^+_\ep v(x_0) \geq v(z) - v(x_0) > u(z) - u(x_0) = \ep S^+_\ep u(x_0).
\end{equation}
Combining \EQ{finite-max-point} and \EQ{temp}, we obtain
\begin{equation*}
S^-_\ep v(x_0) - S^+_\ep v(x_0) < S^-_\ep u(x_0) - S^+_\ep u(x_0).
\end{equation*}
This contradicts \EQ{finite-hyp}, and the lemma follows.
\end{proof}

The next lemma allows us to modify solutions of the PDE \EQ{infinity-laplace} to obtain solutions of the finite difference inequalities \EQ{finite-hyp}. We use the notation
\begin{equation*}
u^\ep(x) := \max_{\bar B(x,\ep)} u \mmbox{and} u_\ep(x) := \min_{\bar B(x,\ep)} u, \quad x\in \Omega_\ep,
\end{equation*}
which allows us to write $\ep S^+_\ep u(x) = u^\ep(x)-u(x)$ and $\ep S^-_\ep u(x) = u(x) - u_\ep(x)$.

\begin{lem}
\label{lem:ballmax}
If $u$ is infinity subharmonic in $\Omega$, then
\begin{equation}
\label{eq:ballmaxsub}
S^-_\ep u^\ep(x) - S^+_\ep u^\ep(x) \leq 0 \quad \mbox{for every} \ \ x \in \Omega_{2 \ep},
\end{equation}
 and if $v$ is infinity superharmonic in $\Omega$, then
\begin{equation}
\label{eq:ballmaxsuper}
S^-_\ep v_\ep(x) - S^+_\ep v_\ep(x) \geq 0 \quad \mbox{for every} \ \ x \in \Omega_{2 \ep}.
\end{equation}
\end{lem}

\begin{proof}
We first prove \EQ{ballmaxsub}.  Fix a point $x_0 \in \Omega_{2 \ep}$. Select $y_0 \in \bar B(x_0,\ep)$
and $z_0 \in \bar B(x_0, 2 \ep)$ such that
$u(y_0) = u^\ep(x_0)$ and $u(z_0) = u^{2\ep}(x_0)$.  We have
\begin{align*}
\ep (S^-_\ep u^\ep(x_0)-S^+_\ep u^\ep(x_0)) &= 2u^\ep(x_0) -
(u^\ep)^\ep(x_0)-(u^\ep)_\ep(x_0)\\
& \leq 2u^\ep(x_0)-u^{2\ep}(x_0)-u(x_0)\\
&=2u(y_0)-u(z_0)-u(x_0).
\end{align*}
A simple calculation verifies that the inequality
\begin{equation*}
u(w) \leq u(x_0) + \frac{u(z_0) - u(x_0)}{2\ep} |w - x_0|
\end{equation*}
holds for all $w \in \partial \left( B(x_0, 2 \ep) \setminus
\{x_0\}\right)$. Since $u$ enjoys comparisons with cones from above,
the inequality therefore holds for every $w \in B(x_0, 2 \ep) \setminus
\{x_0\}$, and thus for every $w \in \bar B(x_0,2\ep)$.
Substituting $w = y_0$ and using the fact that $|y_0 - x_0| \leq \ep$,
we obtain $2 u(y_0) - u(x_0) - u(z_0)\leq 0$. This completes the
proof of \EQ{ballmaxsub}.

To obtain \EQ{ballmaxsuper}, apply \EQ{ballmaxsub} to $-v$ and simplify, using the fact that $(-v)^\ep=-v_\ep$.
\end{proof}


\begin{jensen}[\cite{Jensen:1993}]
Assume that $u, v \in C(\bar \Omega)$ are infinity subharmonic and superharmonic, respectively. Then
\begin{equation}
\label{eq:comparison}
\max_{\bar \Omega} (u - v) = \max_{\partial \Omega} (u - v).
\end{equation}
\end{jensen}

\begin{proof}
According to Lemmas \ref{lem:finite} and \ref{lem:ballmax},
\begin{equation*}
\sup_{\Omega_\ep} (u^\ep -  v_\ep) = \sup_{\Omega_\ep \setminus \Omega_{2 \ep}} (u^\ep - v_\ep)
\end{equation*}
for every $\ep > 0$. We pass to the limit $\ep \to 0$ to obtain \EQ{comparison}.
\end{proof}

\begin{furtherremarks}
M. Crandall \cite{Crandall:personal} has pointed out that our arguments apply, nearly verbatim, if we replace the Euclidean norm in the definition of cone function with a general norm. In this case, the PDE \EQ{infinity-laplace} is replaced by the notion of an \emph{absolutely minimizing Lipschitz extension} (we refer to \cite{Aronsson:2004} for the definition and background). We thereby obtain a new proof of \cite[Theorem 1.4]{Peres:2009} for bounded domains, which establishes the uniqueness of absolutely minimizing Lipschitz extensions in this general setting. The hypotheses of \LEM{finite} can be further relaxed to metric spaces with mild structural conditions.

Our use of the functions $u^\ep$ and $u_\ep$ appears more natural when we recall that a function $u\in C(\Omega)$ is infinity subharmonic if and only if the map $\ep \to u^\ep(x)$ is convex for every $x\in \Omega$ (see \cite[Lemma 4.1]{Crandall:2008}). Crandall \cite{Crandall:personal} has shown that our proof of Jensen's theorem can be efficiently presented by using this equivalence in lieu of \THM{CEG}.

Generalized versions of Lemmas \ref{lem:finite} and \ref{lem:ballmax} appear in \cite{Armstrong:preprint}, wherein the idea of passing from a PDE to a finite difference inequality by maxing over $\ep$-balls is used to obtain new results for the (1-homogeneous) infinity Laplace equation with a nonzero function $f$ on the right-hand side. The finite difference equation $S^+_\ep u = S^-_\ep u$ has been previously considered on finite graphs by Oberman \cite{Oberman:2005} and Le Gruyer \cite{LeGruyer:2007}, both of whom proved existence and uniqueness of solutions to the discretized finite difference equation (with given boundary data).
\end{furtherremarks}

The authors warmly thank Michael Crandall for his many valuable suggestions and remarks. This short article was also improved by the helpful comments of Stephanie Somersille, Kelli Talaska, and Yifeng Yu.

\bibliographystyle{amsplain}
\bibliography{jensen}
\end{document}